\documentclass[preprint,twocolumn]{autart}

\newcommand{\mc}{\mathcal}

\newcommand{\real}{\mathbb{R}}

\newcommand{\realpos}{\mathbb{R}_{\geq 0}}

\newcommand{\transpose}{\mathsf{T}}  
\newcommand{\tsp}{\mathsf{T}}

\newcommand{\zero}{\boldsymbol{0}}
\newcommand{\Pm}{\mathrm{P}}
\newcommand{\Wm}{\mathrm{W}}
\newcommand*{\QEDB}{\hfill\ensuremath{\square}}%
\newcommand*{\QED}{\hfill\ensuremath{\blacksquare}}

\newcommand{\map}[3]{#1: #2 \rightarrow #3}
\newcommand{\setdef}[2]{\{#1 \; : \; #2\}}
\newcommand{\sbs}[2]{{#1}_{\textup{#2}}}
\newcommand{\sps}[2]{{#1}^{\textup{#2}}}
\newcommand{\subscr}[2]{{#1}_{\textup{#2}}}
\newcommand{\supscr}[2]{{#1}^{\textup{#2}}}

\newcommand{\blue}[1]{{\color{blue} #1}}

\usepackage{hyperref}
\usepackage{graphicx,color}
\graphicspath{{./IMG/}{images/}}
\usepackage{amsmath}
\usepackage{amssymb}
\usepackage{mathrsfs}
\usepackage{subfigure}
\usepackage{url}
\usepackage{booktabs}
\usepackage{array}
\usepackage[table]{xcolor}
\usepackage{mathtools} 		
\usepackage{epstopdf}
\usepackage{cite}

\usepackage{graphicx}
\usepackage{amssymb}

\usepackage{lineno}

\journal{Automatica}

\begin{document}

\begin{frontmatter}




\title{\LARGE Secure Trajectory Planning Against \\
Undetectable Spoofing Attacks\thanksref{footnoteinfo}}

\thanks[footnoteinfo]{This paper was not presented at any IFAC
  meeting. Corresponding author: G.~Bianchin. This material is based
  upon work supported in part by ARO award 71603NSYIP, and in part by
  NSF award CNS-1646641.}

\author{Yin-Chen Liu, Gianluca Bianchin, and Fabio Pasqualetti}

\address{Department of Mechanical Engineering, University of
  California, Riverside,
  \{\href{mailto:yliu138@ucr.edu}{\texttt{yliu}},
  \href{mailto:gianluca@engr.ucr.edu}{\texttt{gianluca}},
  \href{mailto:fabiopas@engr.ucr.edu}{\texttt{fabiopas\}@engr.ucr.edu}}.}

\begin{abstract}
  This paper studies, for the first time, the trajectory
  planning problem in adversarial environments, where the objective is to
  design the trajectory of a robot to reach a desired final state
  despite the unknown and arbitrary action of an attacker. In
  particular, we consider a robot moving in a two-dimensional space
  and equipped with two sensors, namely, a Global Navigation Satellite
  System (GNSS) sensor and a Radio Signal Strength Indicator (RSSI)
  sensor. The attacker can arbitrarily spoof the readings of the GNSS
  sensor and the robot control input so as to maximally deviate its
  trajectory from the nominal precomputed path. We derive explicit and
  constructive conditions for the existence of undetectable attacks,
  through which the attacker deviates the robot
  trajectory in a stealthy way. 
  Conversely, we characterize the
  existence of secure trajectories, which guarantee that the robot
  either moves along the nominal trajectory or that the attack remains
  detectable. We show that secure trajectories can only exist between
  a subset of states, and provide a numerical mechanism to compute
  them. We illustrate our findings through several numerical studies,
  and discuss that our methods are applicable to different models of
  robot dynamics, including unicycles. 
  More generally, our results show how control design
  affects security in systems with nonlinear dynamics.
\end{abstract}


\end{frontmatter}


\section{Introduction}
Autonomous robots have rapidly been adopted in a broad range of
applications, including delivery, exploration, surveillance, and
search and rescue. Autonomous robots rely on sensory data to make
decisions, plan their trajectories, and apply controls. Yet, as
demonstrated by recent studies and real world incidents, sensory data
can be accidentally and maliciously compromised, thus undermining the
effectiveness of autonomous operations in critical and adversarial
applications.

Despite recent advances in understanding and enhancing the security of
cyber-physical systems, security tools for autonomous systems are
still of limited applicability and effectiveness. In this paper we
formulate and solve a \emph{secure trajectory planning} problem, where
the objective is to design the trajectory of a robot to reach a
desired final state despite unknown and arbitrary attacks. We consider
a robot equipped with a Global Navigation Satellite System (GNSS)
sensor and a Radio Signal Strength Indicator (RSSI) sensor, and focus
on attackers capable of simultaneously spoofing the GNSS readings and
sending falsified control inputs to the robot. 
We show how the attacker can generate undetectable attacks that 
maximally deviate the robot from the nominal and precomputed 
path, and study how the trajectory planner can exploit the RSSI readings 
to reveal certain attacks. 
Moreover, we demonstrate the existence of secure trajectories between 
certain initial and final configurations, and propose a technique to determine the corresponding control inputs.
We remark that, because of the nonlinearity of RSSI sensor readings, 
existing security methods based on linear models are inapplicable in our 
setting. 
In fact, our results show for the first time that the security of a system with 
nonlinear dynamics can be improved by appropriately designing its control 
inputs.

\noindent
\textbf{Related work} Security of cyber-physical systems is, by now, a
widely studied topic across the controls and computer science
communities, among others. Yet, most methods are applicable to static
systems or systems with linear dynamics
\cite{YZL-AD-FS-IM-MDDB:19,FP-FD-FB:10y,CZB-FP-VG:16,
  YM-BS:10a,FH-PT-SD:11}, and theoretical results and tools for the
security of systems with nonlinear dynamics are still critically
lacking. Few exceptions are \cite{JPH-SDB:19}, which considers
  the problem of controlling a system under attack in a game-theoretic
  framework, \cite{YS-PN-NB-AS-SA-PT-15}, which focuses on nonlinear
  systems satisfying differential flatness properties, and the recent
  articles \cite{QU-DF-YHC-CJT:18,JK-CL-HS-YE-JHS:18}, which are
  however restricted to state the estimation problem in the presence of
  attacks modifying the system measurements only.  Instead, in this
work we focus on characterizing detectability of attacks modifying
both the measurements and the input of the system, on quantifying their 
effects on the trajectories, and on the problem of designing nominal
control inputs to restrict or prevent undetectable attacks.

The literature on GNSS spoofing attack mechanisms and their detection is 
also related to this paper. 
Existing approaches to identify and prevent spoofing attacks
can be divided into two categories: filtering-based and
redundancy-based techniques. 
Filtering-based detection techniques rely on signal processing 
methods to reveal compromised streams of sensory data (e.g., see 
\cite{AB-AJ-VD-JN-GL:12,XJ-JZ-BJH-JJM-ADD:13}). 
Redundancy-based techniques, instead, rely on the availability of 
measurement from multiple types of sensors to reveal inconsistency in the 
data (e.g., see
\cite{PFS-BNA-KCS-RJH:14,QZ-SH-FL-MC:16,DSR-PFS-KCS:15,
LMP-BWO-SPP-JAB-KDW-TEH:14, 
MLP-BWO-JAB-DPS-TEH:13,PYM-TEH-NML:09}).
The methods developed in this work combine these two principles.
In fact, detection is achieved by processing the sensory data over time, 
thus ensuring compatibility between the measurements and the robot 
dynamical model, and by processing the measurements of two or more 
sensors, thus exploiting redundancy across the two channels.

\noindent
{\textbf{Paper contribution} This paper features four main
contributions. First, we demonstrate the existence and characterize the form 
of undetectable attacks, that is, coordinated attack inputs that deviate the 
robot trajectory from the nominal path and cannot be detected using the 
GNSS and RSSI readings. 
Second, we demonstrate how an attacker can design optimal undetectable 
attacks that maximally deviate the robot from its nominal path while 
maintaining undetectability.
Third, we show the existence of secure trajectories where, independently 
of the intensity of the attack, the robot either follows the nominal 
precomputed path or readily detects the attack. 
Fourth, we formulate and solve the secure trajectory planning problem, 
which asks for the design of open-loop control inputs that allow the robot 
to securely navigate from a given initial configuration to a certain final 
position.
As a minor contribution, we study and characterize undetectable attacks 
and secure trajectories for robots with unicycle dynamics, thus showing
that our techniques are in fact applicable to different nonlinear dynamical 
robot models and sensors.
More generally, our results show that secure trajectories can be
substantially different from minimum-time trajectories, and thus 
demonstrate that the security of a system with nonlinear dynamics  
depends upon the inputs to the system adopted for control.

\noindent
\textbf{Paper organization} The remainder of the paper is organized as
follows. Section \ref{sec:problemSetup} presents our problem setup and
attack model. Section \ref{sec:undetectableAttacks} contains our
notion of undetectability and our characterization of undetectable
attacks. Section \ref{sec:attacker} and Section
\ref{sec:SecureTrajectories} contain, respectively, the design of
optimal undetectable attacks and of secure trajectories. Finally,
Section \ref{sec: unicycle} contains an extension of our results to the
case of robots with unicycle dynamics, and Section
\ref{sec:conclusions} concludes the paper.

\section{Problem setup and preliminary
  notions}\label{sec:problemSetup}
We consider a robot with double-integrator dynamics,
\begin{align}\label{eq: nominal dynamics}
\underbrace{
\begin{bmatrix}
\subscr{\dot p}{n}\\
\subscr{\dot v}{n}
\end{bmatrix} 
}_{\subscr{\dot x}{n}}
= 
\underbrace{
\begin{bmatrix}
\zero_2 & I_2\\
\zero_2 & \zero_2
\end{bmatrix}
}_{A}
\underbrace{
\begin{bmatrix}
\subscr{p}{n}\\
\subscr{v}{n}
\end{bmatrix} 
}_{\subscr{x}{n}}
+
\underbrace{
\begin{bmatrix}
\zero_2\\
I_2
\end{bmatrix}  
}_{B}
\subscr{u}{n}
,
\end{align}
where $\map{\subscr{p}{n}}{\realpos}{\real^2}$ denotes the robot
position, $\map{\subscr{v}{n}}{\realpos}{\real^2}$ the robot velocity,
and $\map{\subscr{u}{n}}{\realpos}{\real^2}$ the nominal control input
that actuates the acceleration of the robot. The control input
$\subscr{u}{n}$ is the design parameter that is used to plan the
nominal robot trajectory between two desired configurations. We assume
that $\subscr{u}{n}$ is piecewise continuous, and that
\begin{align*}
\| u_\text{n}\| \leq \subscr{u}{max},
\end{align*}
where $\subscr{u}{max} \in \real_{>0}$.
%
%
We let the robot be equipped with two noiseless sensors: a GNSS
receiver that provides an absolute measure of the position, and a RSSI
sensor that provides a measure of the relative distance between the
robot and a base station located at the origin of the reference
frame. Specifically, the sensor readings are
\begin{align}\label{eq:nominal_readings}
\supscr{y}{GNSS}_\text{n}  = \subscr{p}{n} , \quad \text{and}\quad
\supscr{y}{RSSI}_\text{n}  = \subscr{p}{n}^\tsp \subscr{p}{n} .
\end{align}
Although our results can be extended to include different sensors, we
focus on GNSS and RSSI sensors because they are available in many
practical applications \cite{MJ-RD:03}.

We consider scenarios where the robot operates in an adversarial
environment, where adversaries can simultaneously (i) spoof the
GNSS signal $\supscr{y}{GNSS}_\text{n}$, and (ii) compromise the
control input $u_\text{n}$. The robot dynamics in the presence of
attacks are
\begin{align}\label{eq: dynamics with attack}
  \underbrace{
  \begin{bmatrix}
    \dot{p}\\
    \dot{v}
  \end{bmatrix} 
  }_{\dot x}
  = 
  \underbrace{
  \begin{bmatrix}
    \zero_2 & I_2\\
    \zero_2 & \zero_2
  \end{bmatrix}
              }_A
              \underbrace{
              \begin{bmatrix}
                p\\
                v
              \end{bmatrix} 
  }_x
  +
  \underbrace{
  \begin{bmatrix}
    \zero_2\\
    I_2
  \end{bmatrix}  
  }_B
  u,
\end{align}
where $u \in \real^2$ denotes the attacked control input that obeys
the bound on the maximum acceleration $\|u\| \leq u_\text{max}$, and
\begin{align}\label{eq:attack_readings}
y^\text{GNSS} = p + \supscr{u}{GNSS},\quad \text{and} \quad
y^\text{RSSI} = p^\tsp p ,
\end{align}
where $\map{\supscr{u}{GNSS}}{\realpos}{\real^2}$ denotes the GNSS
spoofing signal.
Finally, we make the practical assumption
that, at time $t=0$, the nominal and attacked configurations satisfy
$\subscr{x}{n}(0) = x(0)$.

In the remainder of this paper, we will denote the state by
$x = [p^\tsp, v^\tsp]^\tsp$ or $p$ and $v$ interchangeably, depending on 
the context. 
In particular, we let
$x = [x_1, x_2, x_3, x_4]^\tsp$, and 
$p = [p^x, p^y]^\tsp = [x_1, x_2]^\tsp $,
$v = [v^x, v^y]^\tsp = [x_3, x_4]^\tsp $.

\begin{rem}{\bf \textit{(Spoofing attack mechanism)}}
Well-known vulnerabilities of GNSS are conventionally associated with 
the lack of appropriate encryption in the signals that are broadcast
by the satellite system.
A typical framework to cast spoofing attacks 
consists in a receiver-spoofer antenna \cite{DPS-TEH-AAF:12} that 
is capable of sensing the authentic GNSS signals and of rebroadcasting 
falsified streams of information at a higher signal intensity.
The retransmitted signals are typically designed in a way to induce the 
GNSS  receiver to resynchronize with the spoofed information, for instance 
by gradually increasing the intensity of the retransmission.
Once the onboard receiver has resynchronized with the falsified signals, 
an attacker can arbitrarily decide the GNSS data received
by the robot, resulting in \eqref{eq:attack_readings}.
A more in-depth discussion of common spoofing schemes and the 
required hardware can be found in e.g. 
\cite{DPS-TEH-AAF:12,AJK-DPS-JAB-TEH:14}.

In common mobile robotic applications, robots communicate wirelessly 
with a ground control station, that is responsible to compute the actions 
and control inputs to be executed by the robot.
The use of wireless communication constitutes a possible vulnerability 
that can be modeled as in \eqref{eq: dynamics with attack}. 
See \cite{KH-CS-13} for a discussion of common vulnerabilities of wireless 
communication in commercial Unmanned Aerial Vehicles (UAV).
\QEDB
\end{rem}

In this work we consider two problems that formalize the contrasting
objectives of an attacker, that is, to deviate the robot trajectory
while remaining undetected, and the trajectory planner, that is, to
design a trajectory between two configurations that is robust to
attacks. We refer to the latter problem to as the 
\emph{secure trajectory planning} problem. 
In particular, the attacker aims to design the attack inputs 
$(u, \supscr{u}{GNSS})$ so that
\begin{enumerate}
\item[(i)] the deviation between the robot nominal trajectory and the actual 
  (attacked) trajectory is maximized; and
\item[(ii)] the attack remains undetected (as defined below).
\end{enumerate}
Instead, the secure trajectory planning problem asks for a
nominal control input $\subscr{u}{n}$ to guarantee that
\begin{enumerate}
\item[(iii)] in the absence of attacks, $\subscr{u}{n}$ allows the robot to
 reach a desired final state; and
\item[(iv)] in the presence of attacks, attacks are detectable (see
  below) by processing the signals $\subscr{u}{n}$,
  $\supscr{y}{GNSS}$, and $\supscr{y}{RSSI}$.
\end{enumerate}
Two observations are in order.
First, although the problem of trajectory planning has a long history in 
robotics (e.g., see \cite{SML:06}), the problem of designing
trajectories in adversarial environments has not been studied before.
Second, the large body of literature on detection and mitigation of
attacks in cyber-physical systems with linear dynamics (e.g., see
\cite{FP-FD-FB:10y,CZB-FP-VG:16}) is not applicable to the considered
secure trajectory planning problem, since the system model is nonlinear due 
to \eqref{eq:attack_readings}. As we show later in this paper, 
and differently from the case of systems with linear dynamics, attack 
detectability for nonlinear systems depends also on the control 
input adopted by the trajectory planner.

We next formalize the notion of attack undetectability.
\begin{defn}{(\bf \textit{Undetectable
      attack})}\label{def: undetectable through y1 and y2}
  The attack $(u, \supscr{u}{GNSS})$ is undetectable if the
  measurements satisfy, at all times,
  \begin{align*}
    \supscr{y}{GNSS} =  \supscr{y}{GNSS}_\text{n}, \;\;\text{ and }\;\;
    \supscr{y}{RSSI} =  \supscr{y}{RSSI}_\text{n} ,
  \end{align*}
  and it is detectable otherwise.
  \QEDB
\end{defn}
Loosely speaking, an attack is undetectable if the measurements
generated by the attacker are compatible with their nominal counterparts
and with the robot dynamics at all times. On the other hand, if the
conditions in Definition \ref{def: undetectable through y1 and y2}
are not satisfied, then the attack is readily detected by simple
comparison between the nominal and actual measurements.

  \begin{rem}{\bf \textit{(Attack detectability in the presence of
        noise)}} In this work we characterize undetectable attacks and
    secure trajectories for deterministic systems. When the dynamics
    or the sensors are driven by noise, different and more relaxed
    notions of attack detectability should be adopted, as done for
    instance in \cite{CZB-FP-VG:16} for the case of linear
    dynamics. Loosely speaking, undetectable attacks are easier to
    cast in stochastic systems, because an attacker has the additional
    possibility of hiding its action within the noise limits. Thus,
    the conditions derived in this paper for deterministic systems
    serve as fundamental limitations also for stochastic systems.
    \QEDB
\end{rem}

Finally, we combine the objectives of the attacker and of the trajectory 
planner into an optimization problem of the general form:
\begin{align}
\label{opt: informal}
& \underset{u,u^\textup{GNSS}}{\text{max }} \;\;\;\;
\underset{\subscr{u}{n},T}{\text{min}} \;\;
& & \int^T_0 L(\subscr{x}{n}, x, t)~dt+V(\subscr{x}{n}(T), x(T)), \nonumber \\
& \text{ \;\;\;subject to:} & & \text{Dynamics } \eqref{eq: nominal dynamics} 
\text{ and } \eqref{eq: dynamics with attack}, \nonumber \\
& & &  (u,\supscr{u}{GNSS})  \text{ is undetectable},
\end{align}
where $T \in \real_{>0}$ represents the planning horizon,
$\map{L}{\real^4 \times \real^4 \times \realpos}{\realpos}$ is an
integral cost, and $\map{V}{\real^4 \times \real^4}{\realpos}$ is a
terminal cost that is chosen to penalize deviations between the
nominal and attacked trajectories at the final time. 
The optimization problem \eqref{opt: informal} captures the
  general class of problems that can be solved with the framework
  proposed in this paper, and will be further specified and discussed
  in the following sections. 
  
We observe that \eqref{opt: informal} is composed of two sequential phases. 
In the first phase, the trajectory planner designs the
  nominal control input $\sbs{u}{n}$ and the control horizon $T$
  (inner minimization problem) to satisfy the objectives (iii) and
  (iv). In the second phase, the attacker designs the attack inputs
  $(u,\sps{u}{GNSS})$ given the nominal input $\sbs{u}{n}$ (outer
  maximization problem) to satisfy objectives (i) and (ii).
Further, we note that  \eqref{opt: informal} can be interpreted as a
  Stackelberg game \cite{TB-GJO:99}, where undetectable attacks 
  represent the best response among all strategies that can be adopted by 
  the attacker, and secure trajectories represent the strategy that maximizes 
  the payoff of the trajectory planner, anticipating the fact that the attacker 
  will adopt its best response.

  \begin{rem}{\bf \textit{(Control mechanism and attacker
        information)}}\label{rem:openLoopControlMechanism}
    Our formulation \eqref{opt: informal} reflects a control framework
    where the trajectory of the robot is planned in an open-loop
    fashion at the beginning of the control horizon by a remote
    control station, and the resulting control parameters are then
    transmitted in batch to the robot \cite{MJ-RD:03}. Our assumptions
    are motivated by the vulnerabilities of wireless communication,
    through which an attacker can intercept the information
    transmitted to the robot. Thus, to successfully cast undetectable
    attacks, the attacker is required to know the robot dynamics and
    the nominal trajectory ahead of time. These requirements can be
    relaxed in the case of single integrator dynamics
    \cite{GB-FP:19}. \QEDB
\end{rem}

\begin{rem}{\bf (\textit{Undetectability with GNSS sensor})}
\label{rem:undetectability through y1}
In scenarios where GNSS is the only sensor for detection, an adversary
can deliberately alter the control input while remaining undetected.
To see this, notice that the effect of any attack $u$ can be canceled
from the GNSS reading by selecting
$\supscr{u}{GNSS}=\subscr{p}{n}-p$. Thus, secure trajectories for the
considered attack model do not exist if the robot has no redundancy 
to combine with the GNSS readings. 
 \QEDB
\end{rem}

\section{Characterization of undetectable
  attacks}\label{sec:undetectableAttacks}
In this section we characterize the class of undetectable attacks and
the resulting attacked trajectories. First, we establish a
relationship between the nominal and attacked instantaneous position
and velocity. Then, we derive an explicit expression of undetectable
attacks, and demonstrate how an attacker can readily design attacks
that escape detection. The following result relates attack
trajectories with their nominal counterparts.

\begin{lem}{\bf \textit{(Undetectable
      trajectories)}}\label{thm:undetectable_trajectories}
  Let $(u,\supscr{u}{GNSS})$ be an undetectable attack.  Then,
  \begin{align*}
    p^\transpose p=\subscr{p}{n}^\transpose\subscr{p}{n}, \text{ and }
    v^\transpose p=\subscr{v}{n}^\transpose \subscr{p}{n} .
  \end{align*}
\end{lem}
\begin{pf}
  The first equality in the statement follows by substitution of
  \eqref{eq:nominal_readings} and \eqref{eq:attack_readings} into
  Definition~\ref{def: undetectable through y1 and y2}.  Further, by
  taking the time-derivative on both sides of the equality
  $p^\transpose p = \subscr{p}{n}^\transpose \subscr{p}{n} $, and by
  using the assumption $\subscr{x}{n}(0) = x(0)$, we obtain
  $2 \dot p^\transpose p = 2\subscr{\dot p}{n}^\transpose \subscr{p}{n} $
  at all times, from which the statement follows.
\QED
\end{pf}

From Lemma~\ref{thm:undetectable_trajectories}, trajectories
generated by undetectable attacks are characterized by two features:
at all times, (i) the distance $p^\transpose p$ between the attacked
position and the RSSI-base station must equal the distance
$\subscr{p}{n}^\transpose \subscr{p}{n}$ in the nominal trajectory,
and (ii) the component of the velocity $v$ along the position $p$ must
equal the component of the nominal velocity $\subscr{v}{n}$ along
$\subscr{p}{n}$. These two geometric properties are illustrated in
Fig.~\ref{fig:conditions_thm1}. Next, we give an implicit
characterization of undetectable attacks.

\begin{figure}[tb]
\centering
\includegraphics[width=.6\columnwidth]{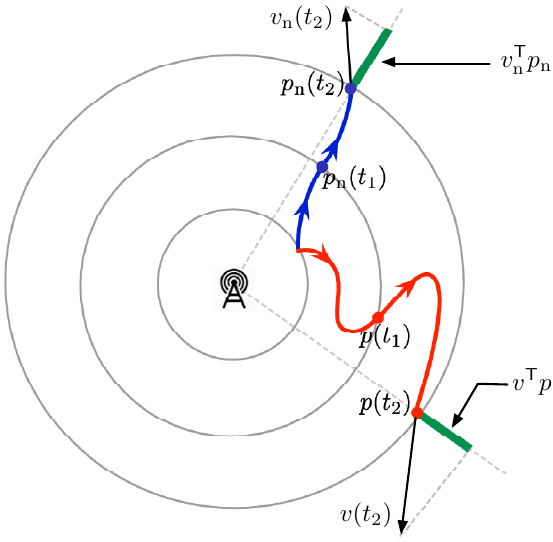}
\caption{Nominal (blue) and undetectable attack (red) trajectories. 
  At all times, the two trajectories have equal distance from the base
  station, and equal velocity components along the direction of
  the instantaneous position (green segments).  }\label{fig:conditions_thm1}
\end{figure}

\begin{thm}{\bf \textit{(Implicit characterization of undetectable
      attacks)}}\label{thm: main zero dynamics}
  The attack $(u,\supscr{u}{GNSS})$ is undetectable if and only if
  \begin{align}\label{eq: implicit undetectability condition}
    u^\tsp p  = \subscr{u}{n}^\tsp \subscr{p}{n} 
    + \|\subscr{v}{n} \|^2  - \| v \|^2, \text{ and }
    \supscr{u}{GNSS} = \subscr{p}{n}-p.
  \end{align}
\end{thm}
\begin{pf}
  \textit{(Only if)} By substitution of \eqref{eq:nominal_readings}
  and \eqref{eq:attack_readings} into Definition~\ref{def:
    undetectable through y1 and y2} we readily obtain
  $p + \supscr{u}{GNSS} = \subscr{p}{n}$, from which the second
  identity in the statement follows.
To show the first identity, we note that for every undetectable attack $u$ the
identity $\supscr{y}{RSSI} - \supscr{y}{RSSI}_\text{n} = 0$ holds.
We then make explicit the dependence on the control input in the above 
identity by taking the second-order derivative with respect to time. This
yields $  \supscr{\ddot y}{RSSI} -  \supscr{\ddot y}{RSSI}_\text{n} =0$.
Then,
\begin{align*}
0 =  \supscr{\ddot y}{RSSI} -  \supscr{\ddot y}{RSSI}_\text{n} = 
2 u^\transpose p + 2 v^\transpose v -
2  \subscr{u}{n}^\transpose \subscr{u}{n} - 2 \subscr{v}{n}^\transpose
    \subscr{v}{n},
\end{align*}
from which \eqref{eq: implicit undetectability condition} follows.
We emphasize that the functions 
$\supscr{\dot y}{RSSI}$, 
$\supscr{\ddot y}{RSSI}$, 
$\supscr{\dot y}{RSSI}_\text{n}$, 
$\supscr{\ddot y}{RSSI}_\text{n}$ 
are piecewise continuous functions, since 
the signals $\supscr{ y}{RSSI}$ and  $\supscr{y}{RSSI}_\text{n}$ 
are continuous and twice differentiable at all times. 
To see this, we combine the piecewise continuity assumption on
$\subscr{u}{n}$ and $u$ with the dynamical equations 
\eqref{eq: nominal dynamics} and \eqref{eq: dynamics with attack}, 
and note that the velocities $\subscr{v}{n}$ and $v$ and the positions 
$\subscr{p}{n}$ and $p$ are continuous and differentiable functions of time. 

\noindent
\textit{(If)}
Let $  \supscr{u}{GNSS} = \subscr{p}{n}-p$. 
Substituting into \eqref{eq:attack_readings} yields
  \begin{align*}
    \supscr{y}{GNSS} - \supscr{y}{GNSS}_\text{n}= 
    p + \supscr{u}{GNSS} - \subscr{p}{n} = 0,
  \end{align*}
from which the first condition in 
Definition~\ref{def: undetectable through y1 and y2} follows.
To prove RSSI undetectability, let $u$ satisfy 
\eqref{eq: implicit undetectability condition}. Then,
\begin{align*}    
\supscr{\ddot y}{RSSI} -  \supscr{\ddot{y}}{RSSI}_\text{n} =
2 u^\transpose p + 2 v^\tsp v -  
  2 \subscr{u}{n}^\tsp \subscr{p}{n} - 2 \subscr{v}{n}^\tsp \subscr{v}{n} =0,
\end{align*}
from which we readily obtain the identity
$\supscr{\ddot y}{RSSI} -  \supscr{\ddot{y}}{RSSI}_\text{n} = 0$.
Since the functions $\supscr{y}{RSSI}$ and $\supscr{y}{RSSI}_\text{n}$
are continuous and twice differentiable, and the initial conditions satisfy 
$p(0) = p_\text{n} (0)$ and $v(0) = v_\text{n} (0)$ we conclude that
$\supscr{\dot y}{RSSI} -  \supscr{\dot{y}}{RSSI}_\text{n} = 0$ and 
$\supscr{y}{RSSI} -  \supscr{y}{RSSI}_\text{n} = 0$. 
Thus, $\supscr{y}{RSSI} = \supscr{y}{RSSI}_\text{n}$ at all times, which 
implies undetectability of the attack $u$ and concludes the proof.
\QED
\end{pf}

Finally, we exploit Theorem~\ref{thm: main zero dynamics} to give an
explicit and comprehensive characterization of undetectable attacks.

\begin{cor}{\bf \textit{(Explicit characterization of
      undetectable attacks)}}\label{cor: undetectable attack}
  The attack $(u,\supscr{u}{GNSS})$ is undetectable if and only if it
  satisfies
  \begin{align}\label{eq:undetectable attacker input}
    u = \subscr{a}{r } p + w \text{ and }
                                 \supscr{u}{GNSS} = \subscr{p}{n} - p,
  \end{align}
  whenever $\|p \| \neq 0$, where $w^\tsp p = 0$ and
  \begin{align*}
    \subscr{a}{r} = \frac{\subscr{u}{n}^\tsp \subscr{p}{n} 
    + \|\subscr{v}{n} \|^2  - \| v \|^2}{\| p \|^2} .
  \end{align*}
\end{cor}

Corollary~\ref{cor: undetectable attack} provides a systematic way to
design undetectable attacks by designing the attack inputs 
$(u,\sps{u}{GNSS})$. 
We also note that the input $w$ can be arbitrarily selected by
the attacker and it does not affect detectability of the
attack. Finally, it should be noticed that $\subscr{a}{r}$ corresponds
to the radial acceleration of the robot, that is, the projection of $u$
along $p$, and that the attack input $u$ is unconstrained when
$\Vert p \Vert = 0$ (see also Theorem~\ref{thm: main zero dynamics}).

\section{Design of optimal undetectable attacks}\label{sec:attacker}
In this section we design undetectable attacks that introduce maximal
deviations between the nominal and attacked trajectories. 
Assuming that the attacker knows the nominal control input, we address the
optimal control problem
\begin{subequations}
\label{opt:attacker}
\begin{align}
&\underset{w}{\text{max}}  
					& & \Vert p(T)-p_\text{n}(T) \Vert, \nonumber \\ 
 & \text{subject to} & &  \dot x =Ax+Bu, \label{opt:attacker-a}\\
&&& u = \subscr{a}{r } p + w,  \label{opt:attacker-c}\\
&&& \subscr{a}{r} = (\subscr{u}{n}^\tsp \subscr{p}{n} 
	+ \|\subscr{v}{n} \|^2  - \| v \|^2)\| p \|^{-2},  \label{opt:attacker-d}\\
&&& \|u\|\leq \subscr{u}{max}.
\end{align}
\end{subequations}
In the maximization problem \eqref{opt:attacker}, constraint
\eqref{opt:attacker-a} corresponds to the attacked dynamics \eqref{eq:
  dynamics with attack}, while
\eqref{opt:attacker-c}-\eqref{opt:attacker-d} enforce
attack-undetectability from Corollary~\ref{cor: undetectable attack}. 
We next characterize the optimality conditions
of the problem~\eqref{opt:attacker}.
Let $e_i$ denote the $i$-th canonical vector of appropriate dimension, and 
let $\operatorname{sgn(~)}$ denote the sign function.

\begin{thm}{\bf \textit{(Attack optimality
      conditions)}}\label{thm: optimized attacker input}
  Let $\subscr{a}{r}$ be as in \eqref{opt:attacker-d},
  and let $w^*$ be an optimal solution to \eqref{opt:attacker}.  Then,
  \begin{align*}
    w^* 	= \subscr{a}{t} W x,
  \end{align*}
  where
  \begin{align*}
    \subscr{a}{t} &= - \operatorname{sgn}(\lambda^\tsp B \Wm x)
    \sqrt{\subscr{u}{max}^2 / \| p \|^2 - \subscr{a}{r}^2},\\
    W &= 
        \begin{bmatrix}
          -e_2 & e_1 & 0_2 & 0_2
        \end{bmatrix}
                           ,
  \end{align*}
  and $\lambda$ and $x$ satisfy
\begin{align*}
\dot x &= Ax + \subscr{a}{r} Bp + Bw^*,\\
-\dot \lambda &=      (A^\tsp \lambda + \subscr{a}{r} \tilde{P}
     								+ \subscr{a}{t} \tilde{W} ) \lambda+
     (x^\tsp \tilde{P} \lambda + 2 \subscr{a}{r} \nu \Vert \subscr{p}{n} \Vert^2) 
    				\nabla_x \subscr{a}{r},
\end{align*}
with boundary conditions
\begin{align*}
  x(0) = \subscr{x}{n}(0), \text{ and }
  \lambda(T) = -2[(p(T) - \subscr{p}{n}(T))^\tsp ~ 0_2^\tsp]^\tsp,
\end{align*}
where $\tilde P = P^\tsp B^\tsp$, $\tilde W = W^\tsp B^\tsp$, 
$\nabla_x \subscr{a}{r} = 2 \| p \|^{-2} [0_2^\tsp ~ v^\tsp]^\tsp,$
and 
  $ \nu = - \frac{\lambda^\tsp B \Wm x }{ 2 \subscr{a}{t} \| p \|^2 }$.
\end{thm}

\begin{pf}
To formalize the result, we make use of the fact that any undetectable 
attack \eqref{eq:undetectable attacker input} can be written in the form
\begin{align}
\label{eq:uMatricesML}
u = \subscr{a}{r } \Pm x + \subscr{a}{t } \Wm x,
\end{align}
where $\Pm = [e_1~e_2~0~0] \in \real^{2 \times 4}$, 
$\Wm = [-e_2~e_1~0~0]\in \real^{2 \times 4}$, and 
$\map{\subscr{a}{t}}{\realpos}{\real}$.
Following expression \eqref{eq:uMatricesML}, the function $\subscr{a}{t}$ 
represents the new design parameter in the optimization problem.
To derive the optimality conditions for \eqref{opt:attacker}, we use the 
Pontryagin's Maximum Principle \cite{IMG-RAS:00}, combined with 
the direct adjoining method for mixed state-input constraints 
\cite{RFH-SPS-RGV:95}.
We incorporate \eqref{opt:attacker-c} and  \eqref{eq:uMatricesML} into 
\eqref{opt:attacker-a} and define the \textit{Hamiltonian} 
\begin{align*}
\mc H(x, \subscr{a}{t}, \lambda, t) =  
\lambda^\tsp (Ax + B(\subscr{a}{r} \Pm x + \subscr{a}{t} \Wm x)),
\end{align*}
where 
$\map{\lambda}{[0,T]}{\real^4}$ is a vector function of system costates, 
with the additional constraints 
\begin{align*}
x(0) =\subscr{x}{n}(0), 
\subscr{a}{r} = \frac{\subscr{u}{n}^\tsp \subscr{p}{n} 
	+ \|\subscr{v}{n} \|^2  - \| v \|^2}{\| p \|^2},  
\|u\| \leq \subscr{u}{max}.
\end{align*}
We then use \eqref{eq:uMatricesML}  to rewrite the bound 
$\|u\| \leq \subscr{u}{max}$ as
\begin{align*}
\subscr{a}{r}^2 \| \subscr{p}{n} \|^2   + \subscr{a}{t}^2 \| \subscr{p}{n} \|^2   
\leq \subscr{u}{max}^2,
\end{align*}
and form the Lagrangian by adjoining the Hamiltonian with the 
considered state constraint:
\begin{align*}
\mc L(x, \subscr{a}{t}, \lambda, t, \nu) =  
	\mc H(x, \subscr{a}{t}, &\lambda, t) + \\
	&\nu (\subscr{a}{r}^2 \| \subscr{p}{n} \|^2 + 
	\subscr{a}{t}^2  \| \subscr{p}{n} \|^2- \subscr{u}{max}^2 ),
\end{align*}
where $\map{\nu}{[0,T]}{\real}$ is the Lagrange multiplier associated with 
the state constraint.

By application of the Maximum Principle \cite{RFH-SPS-RGV:95}, the 
optimal control input $\subscr{a}{t}^*$ minimizes the Hamiltonian
over the set
$U(x) = \{\subscr{a}{t} : 
\subscr{a}{r}^2 \|p\|^2 + \subscr{a}{t}^2  \|p\|^2 \leq \subscr{u}{max}^2 \}$, 
that is,
$ \subscr{a}{t}^*=\operatorname{arg} \min_{\subscr{a}{t} \in U(x)} 
			\mc H(x, \subscr{a}{t}, \lambda, t)$.
This fact yields the optimal control law
\begin{align*}
\subscr{a}{t}^* =  -\operatorname{sgn}(\lambda^\tsp B \Wm x) 
\sqrt{\subscr{u}{max}^2 / \|p\|^2 -\subscr{a}{r}^2}.
\end{align*}
Moreover, it follows from the Maximum Principle that there exists  a vector 
function of system costates  $\lambda$ that satisfies the following system 
of equations
\begin{align*}
\dot x &= \frac{\partial \mc L}{\partial \lambda }
\Rightarrow \;\;\; \dot x = Ax + B( \subscr{a}{r } \Pm x +\subscr{a}{t} \Wm x),\\
-  \dot \lambda &= \frac{\partial \mc L}{\partial x } \Rightarrow 
\begin{aligned}[t]
-\dot \lambda = &A^\tsp \lambda + \subscr{a}{r} P^\tsp B^\tsp \lambda +
x^\tsp P^\tsp B^\tsp \lambda \nabla_x \subscr{a}{r} + \\ 
&\subscr{a}{t} W^\tsp B^\tsp \lambda
+ 2 \subscr{a}{r} \nu  \Vert \subscr{p}{n} \Vert^2 \nabla_x \subscr{a}{r}
\end{aligned}\\
0 &= \frac{\partial \mc L}{\partial \subscr{a}{t}} \Rightarrow
					\;\;\; 0 = \lambda^\tsp B \Wm x + 2 \nu \subscr{a}{t} \| p \|^2,
\end{align*}
where 
$\nabla_x \subscr{a}{r}$ denotes the gradient of $\subscr{a}{r}$ with 
respect to $x$, and with boundary conditions 
\begin{align*}
x(0) \;= \subscr{x}{n}(0), \text{ and }
\lambda(T) = \frac{\partial}{\partial x}\Vert p(T)-p_\text{n}(T) \Vert.
\end{align*}
The statement of the theorem follows by substituting the expression of the 
gradient $\nabla_x \subscr{a}{r}  =  2 \| p \|^{-2} [0~0~ v^\tsp]^\tsp$.
\QED
\end{pf}

From Theorem \ref{thm: optimized attacker input}, optimal undetectable
attacks can be computed by solving a two-point boundary value
problem \cite{HBK:18}. 
This class of problems is typically solved numerically, and
it may lead to numerical difficulties for general cases 
\cite{HBK:18}.
To conclude
this section and provide some intuition in the design of optimal
undetectable attacks, we next present an example where optimal attacks
can be characterized explicitly.

\begin{figure}[t]
  \centering \subfigure[]{
     \raisebox{4mm}{
    \includegraphics[width=.44\columnwidth]{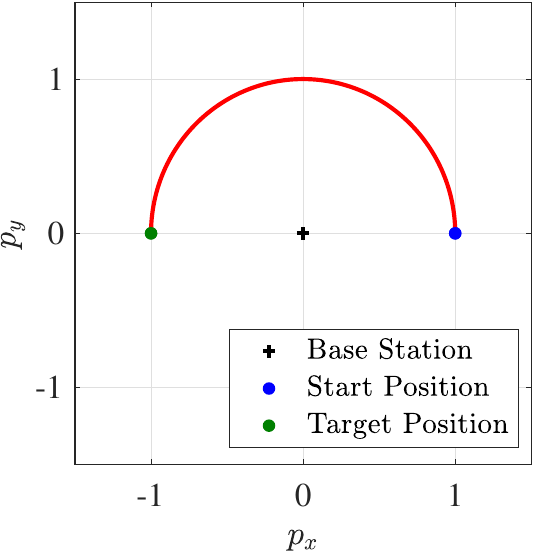} 
    \label{fig:circularMotion_a}  }} 
\begin{minipage}[b]{0.5\columnwidth}
 \subfigure[]{
    \includegraphics[width=.9\columnwidth]{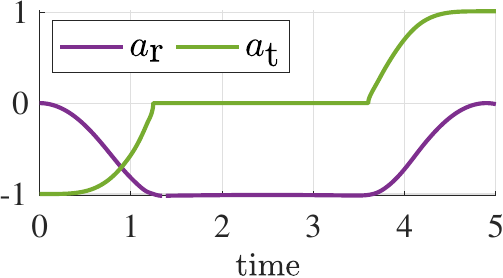} }
\subfigure[]{
    \includegraphics[width=.9\columnwidth]{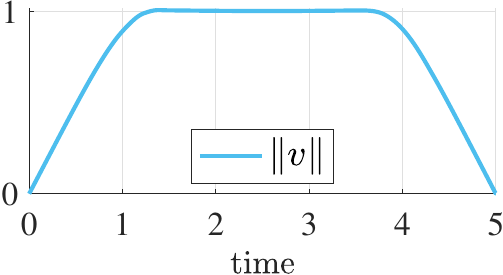} }
\end{minipage}
\caption[]{For an idle robot at position $(1,0)$ (blue dot), (a) shows
  an optimal undetectable attack trajectory, which maintains a
  distance equal to $1$ from the base station and maximises the
  distance from the nominal position. Fig. (b) shows the radial and
  tangential components of the acceleration. Fig. (c) shows the
  velocity of the attacked robot, which becomes zero when the robot
  reaches the final point $(-1,0)$.}
\label{fig:circularMotion}
\end{figure}

\begin{exmp}{\bf \textit{(Undetectable trajectories for idle
      robots)}}
  Let $\subscr{p}{n}(0) = [1~0]^\tsp$, $\subscr{v}{n}(0) = 0$, and
  $\subscr{u}{n} = 0$, so that the robot remains at position
  $\subscr{p}{n}(0)$ at all times. 
Let $\subscr{u}{max} =1$ and $T=5$.
Under these assumptions, the following control inputs satisfy the optimality 
conditions in 
Theorem~\ref{thm: optimized attacker input}:
  \begin{align}
\label{eq:optimalityIdleRobot}
    \subscr{a}{t}^* 	=  
    \begin{cases}
\zeta   \sqrt{\subscr{u}{max}^2  -  \subscr{a}{r}^2},    & t \in [0,\tau], \\
-\zeta    \sqrt{\subscr{u}{max}^2  -  \subscr{a}{r}^2},    & t \in [\tau,T],
    \end{cases} \text{ and }
    \subscr{a}{r} =  - \| v \|^2,
  \end{align}
where $\tau = 3.475$, and $\zeta \in \{-1,1\}$.
It is worth noting that, because at every time the radial acceleration is 
proportional to the square of the magnitude of the velocity, the control 
input~\eqref{eq:optimalityIdleRobot} leads the attacked robot to 
perform a circular motion around the origin (see 
Fig.~\ref{fig:circularMotion}).
Notice that $\zeta = -1$ and $\zeta = 1$ achieve counterclockwise and 
clockwise motion, respectively.
Finally, \eqref{eq:optimalityIdleRobot} is an optimal solution to the 
optimization problem \eqref{opt:attacker} since the deviation 
$\Vert p(T)-p_\text{n}(T) \Vert$ is 
maximized, as graphically illustrated in Fig.~\ref{fig:circularMotion}.

To derive the value of the final time $\tau$,we can explicitly derive a
an expression for the magnitude of the velocity vector $n := \Vert v \Vert$, 
that reads $\dot n = \sqrt{\subscr{u}{max}^2 - n^2}$, where we have 
substituted the expression for $\subscr{a}{r}$ into $\subscr{a}{t}^*$, and 
used the fact that $n$ is independent of $\subscr{a}{r}$.
To obtain the value of $\tau$, we seek for the time needed to steer $n(t)$ 
from $\subscr{u}{max}$ to a full stop, by letting $n(0)= \subscr{u}{max}$
and $n(\tau)= 0$.
\QEDB
\end{exmp}

We show in Fig.~\ref{fig:varyUmax} a set of simulations that illustrate
the effects of optimal attacks (red) when the nominal trajectory is the
shortest path between the initial and the final position (blue).  In
particular, the figure demonstrates that increasing levels of deviation
are achieved by the attacker when the trajectory planner employs control 
inputs with decreasing magnitude (i.e., different and decreasing fractions of 
$\subscr{u}{max}$).

\begin{figure}[tb]
  \centering
  \includegraphics[width=.8\columnwidth]{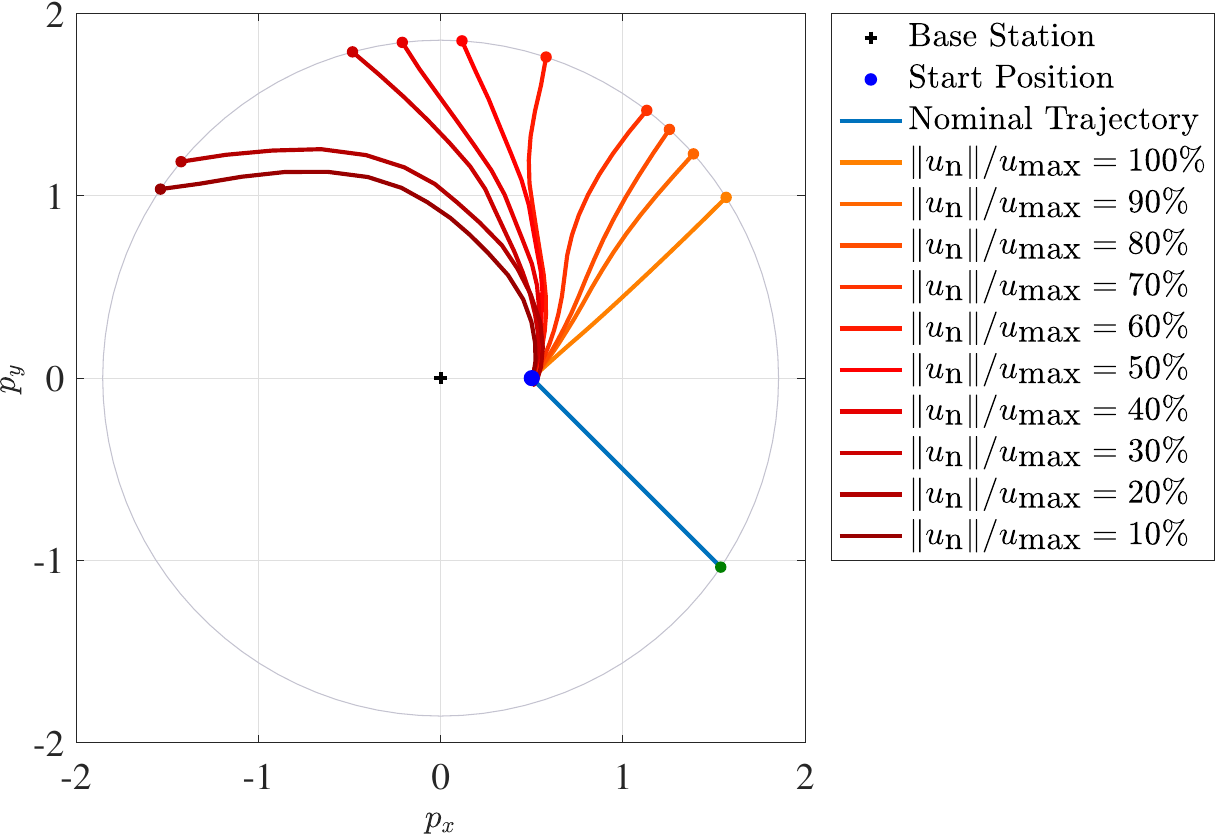}
  \caption{For a nominal straight line trajectory (blue), the figure
    shows the undetectable attack trajectories obtained from Theorem
    \ref{thm: optimized attacker input} for different values of the
    nominal acceleration
    $\Vert \subscr{u}{n} \Vert / \subscr{u}{max}$. As the nominal
    acceleration decreases, the deviation induced by an optimal attack
    increases. Simulation parameters: $T=1.5$, $v(0) = 0$, and
    $\subscr{u}{n} = [1,-1]^\tsp$ at all times.}
  \label{fig:varyUmax}
\end{figure}

\section{Design of secure trajectories}\label{sec:SecureTrajectories}
In this section we address the secure trajectory planning
problem. First, we characterize the existence of secure trajectories
as a function of the initial and final configurations of the
robot. Then, we formulate and solve an optimization problem to design
control inputs that generate secure trajectories to reach a desired final
configuration.
We start with some necessary definitions.
We say that a trajectory $\subscr{x}{n}$ is \emph{secure} if, independently 
of the attack $u$, one of the following mutually exclusive conditions is 
satisfied:
\begin{enumerate}
\item[(C1)] $p=p_\text{n}$ at all times; or\label{cod: security1}
\item[(C2)] if $p\neq p_\text{n}$ at some time, the attack $u$ is
  detectable.\label{cod: security2}
\end{enumerate}
Similarly, a control input is \emph{secure} if it generates a secure
trajectory. We next characterize secure control inputs explicitly.

\begin{thm}{\bf \textit{(Secure control
      inputs)}}\label{thm:secureControlInput}
  Let $\subscr{x}{n}$ be the trajectory generated by
  $\subscr{u}{n}$. Then, $\subscr{u}{n}$ is secure if and only if the 
  following conditions hold simultaneously:
  \begin{enumerate}
\item[(1)] there exists a function $\map{\kappa}{\realpos}{ \{-1,1\} }$ 
	satisfying
    \begin{align}\label{eq:secureControlInput}
      u_\text{n} = \kappa \frac{p_\text{n}}{\|p_\text{n}\|}
      \subscr{u}{max} ,
    \end{align}

  \item[(2)] the trajectory $\subscr{x}{n}$ satisfies $p_\text{n} \neq 0$
    at all times.
  \end{enumerate}
\end{thm}

\begin{pf}
(\textit{If}) 
We assume (1)-(2) and show that either (C1) or 
(C2) is satisfied. 
We distinguish among two cases. \\
\textit{(Case 1)} The attack  $(u,\supscr{u}{GNSS})$ 
does not satisfy the undetectability condition 
\eqref{eq: implicit undetectability condition}.
Then, (C2)  follows. \\
\textit{(Case 2)} The attack $(u,\supscr{u}{GNSS})$ satisfies the 
undetectability condition \eqref{eq: implicit undetectability condition}. 
Under this assumption we now show that (1)-(2) imply (C1).
We first consider the time instant $\tau=0$.
By using the assumption $\subscr{x}{n}(\tau) = x(\tau)$, which yields
$\subscr{p}{n}(\tau) = p(\tau)$ and 
$\Vert \subscr{v}{n}(\tau) \Vert = \Vert v(\tau) \Vert$,
and by substituting into \eqref{eq: implicit undetectability condition} we 
obtain the following undetectability condition valid at time $\tau$:
\begin{align}
\label{eq:undetctabilityCase2}
\subscr{u}{n}^\tsp(\tau) \subscr{p}{n} (\tau) = u^\tsp (\tau) p (\tau).
\end{align}
By taking the $2$-norm on both sides of the above equality, and by 
substituting the expression \eqref{eq:secureControlInput} we obtain
\begin{align*}
\subscr{u}{max} \Vert \subscr{p}{n} (\tau) \Vert 						= 
\vert \subscr{u}{n}(\tau)^\tsp &\subscr{p}{n}(\tau) \vert 			=\\
&\vert u^\tsp(\tau) p(\tau) \vert 												\leq 
\subscr{u}{max} \Vert p (\tau) \Vert,
\end{align*}
where we used the Cauchy-Schwarz inequality. 
Since exact equality must hold, the vectors $u(\tau)$ and $p(\tau)$ are
linearly dependent and $\Vert u (\tau) \Vert = \subscr{u}{max}$, that is,
\begin{align*}
u(\tau) = \gamma(\tau) \frac{p(\tau)}{\Vert p(\tau) \Vert} \subscr{u}{max},
\end{align*}
where $\gamma (\tau ) \in \{-1,1\}$.
Finally, we note that $\gamma(\tau) \neq \kappa (\tau)$ results in a 
violation of \eqref{eq:undetctabilityCase2}, and therefore
$u(\tau) = \subscr{u}{n}(\tau)$.
As a result, $\subscr{p}{n}(\tau^+) = p(\tau^+)$.
To conclude the proof, we iterate the above reasoning for all 
$\tau \in [0,T)$, from which (C1) follows.\\

\noindent
(\textit{Only if})
We show that (C1)-(C2) imply (1)-(2) or, equivalently, if (1)-(2) do not 
simultaneously hold, then (C1)-(C2) are not satisfied.
We distinguish two cases.\\
\textit{(Case 1)}
Let $\subscr{\bar u}{n}$ be any control input that does not satisfy
\eqref{eq:secureControlInput}.
That is, there exists $\bar t \in [0,T]$, such that 
\begin{align*}
 \subscr{\bar u}{n}(t) = \subscr{u}{n}(t) \text{ for all } t \in[0,\bar t), \quad 
 \text{ and } \quad  \subscr{\bar u}{n}(\bar t) \neq u_\text{n}(\bar t).
\end{align*}
Let $\bar u$ be an attack input, with 
$  \bar u = \subscr{\bar a}{r } p +  \bar w$ as in 
Corollary~\ref{cor: undetectable attack}.
We take the absolute value of \eqref{eq:undetectable attacker input} and
use the relationship $\Vert \subscr{v}{n}(\bar t) \Vert = \Vert v(\bar t) \Vert $
to obtain the inequality
\begin{align*}
\vert \subscr{\bar a}{r} (\bar t) \vert =
\frac{\vert \subscr{\bar u}{n}^\tsp(\bar t) \subscr{p}{n}(\bar t) \vert}{
\Vert p(\bar t) \Vert} <
\subscr{u}{max},
\end{align*}
where strict inequality follows from the assumption
$\subscr{\bar u}{n}(\bar t) \neq u_\text{n}(\bar t)$.
As a result, any vector $\bar w$ that satisfies
$ \Vert \subscr{\bar a}{r}(\tau) p(\tau) + \bar w(\tau) \Vert =
\subscr{u}{max}$, 
is a nonzero undetectable attack that violates (C1) and (C2).\\
%
\textit{(Case 2)}
There exists $\bar t \in [0,T]$, such that  $\| \subscr{p}{n}( \bar t) \| = 0$.
It follows from \eqref{eq: implicit undetectability condition}
that whenever $\subscr{p}{n}(\bar t) = 0$ any attack input is unconstrained 
at time instant  $\bar t$. 
As a result, any $u$ with $u(t) = \subscr{u}{n} (t)$ for all 
$t \in [0,\bar t)$ and $u(\bar t) \neq \subscr{u}{n} (\bar t)$ is undetectable 
and violates (C1) and (C2).
\QED
\end{pf}

Theorem \ref{thm:secureControlInput} provides an explicit
characterization of secure control inputs, and it shows that any
secure control input has maximum magnitude $\subscr{u}{max}$ at all
times, and its direction is aligned with the direction of the positioning 
vector. We next show that any secure trajectory evolves on an invariant
manifold of the state space that is uniquely defined by the initial
state of the robot.

\begin{lem}{\bf \textit{(Invariant manifold of secure
      trajectories)}}\label{lem:accessibleConfigurations}
  Let $\subscr{x}{n}$ be the trajectory generated by the secure
  control input $\subscr{u}{n}$. Then, $\subscr{x}{n} \in \mc S$ at
  all times, where
  \begin{align*}
    \mc S=\setdef{x}{x_1 x_4 - x_2 x_3 = x_1(0) x_4(0) - x_2(0)
    x_3(0)}. 
  \end{align*}
\end{lem}
\begin{pf}
  Let $\subscr{x}{n}=[x_1~x_2~x_3~x_4]^\tsp$ denote the solution to
  \eqref{eq: nominal dynamics} with initial condition
  $\subscr{x}{n}(0)$, subject to control inputs that satisfy (1)-(2) in 
  Theorem~\ref{thm:secureControlInput}.
To prove that
  $\subscr{x}{n} \in \mc S$, we equivalently show that the quantity
  $x_1 x_4 - x_2 x_3 $ is time-invariant, that is,
  $\frac{d}{dt}(x_1 x_4 - x_2 x_3) = 0$.  In fact,
  \begin{align*}
    \dot x_1 x_4 +  x_1 \dot x_4 - \dot x_2 x_3 - x_2 \dot x_3=0,
  \end{align*}
  where the last equality follows by substitution of \eqref{eq:
    nominal dynamics}.
\QED
\end{pf}

Lemma \ref{lem:accessibleConfigurations} shows that secure
trajectories are constrained to evolve on a manifold that is defined 
by the initial state of the robot, and it implies that only a subset of the state 
space can be reached via secure trajectories. 
These observations are illustrated in the next example.

\begin{exmp}{\bf \textit{(Reachable
      configurations)}}\label{exmp:accessibleSpace}
  Let $\subscr{x}{n}(0) = [\bar x_1~0~\bar x_3~0]^\tsp$, where
  $\bar x_1 \in \real_{>0}$ and $\bar x_3 \in \real$, that is, let the robot at 
time $t=0$ be located on the horizontal axis, with initial velocity in the 
horizontal direction.
From Lemma~\ref{lem:accessibleConfigurations}, every secure trajectory 
satisfied $\subscr{x}{n} \in \mc S$ at all times, with
  \begin{align*}
    \mc{S}= \setdef{x}{x_1 x_4 - x_2 x_3 = 0 }.
  \end{align*}
It should be observed, however, that secure trajectories may
  in fact be constrained on a strict subset of $\mc S$. In fact, by combining
  the system dynamics \eqref{eq: nominal dynamics} with the secure
  control input \eqref{eq:secureControlInput}, we obtain
  \begin{align*}
    x_1>0, \quad x_2=0, \quad x_4=0,
  \end{align*}
  which is a strict subset of $\mc S$.  
  Notice that the above equations imply that the motion of the robot under 
  secure control inputs is constrained on the positive $x$-axis.
  \QEDB
\end{exmp}

To determine a secure trajectory from the initial position $\subscr{p}{I}$ 
with given velocity $\subscr{v}{I}$ towards the final position 
$\subscr{p}{F}$, we consider the optimization problem\footnote{In the 
optimization problem, we allow a free final velocity to ensure that the final 
configuration achieved belongs to the invariant manifold that constraints 
the secure trajectory.}
\begin{align}
  &\underset{\kappa, T}{\text{min}}  
  & & T ~+~ 
      \Vert \subscr{p}{n}(T)-\subscr{p}{F}\Vert, \nonumber \\ 
  & \text{subject to} & & \subscr{\dot x}{n} =A \subscr{x}{n} + 
  		B \subscr{u}{n}, \nonumber\\
  &&& \subscr{p}{n}(0) = \subscr{p}{I}, \nonumber 
  \subscr{v}{n}(0) = \subscr{v}{I},\nonumber \\
  &&& u_\text{n} = \kappa \frac{p_\text{n}}{\|p_\text{n}\|}
      u_\text{max} ,
      \label{opt:secureTrajectoryPlanning}
\end{align}
which aims to find a secure control input that minimizes a weighted
combination of the distance to the desired final position and the time
needed to reach such position.
The following result characterizes the solutions to the
optimization problem~\eqref{opt:secureTrajectoryPlanning}.

\begin{thm}{\bf \textit{(Optimality conditions for secure control
      inputs)}}\label{thm:secureTrajectoryPlanning}
  Let $\kappa^*$ and $T^*$ be an optimal solution to
  \eqref{opt:secureTrajectoryPlanning}. Then,
  \begin{align*}
    \kappa^* = -\operatorname{sgn}(\lambda^\tsp B \subscr{p}{n}), \text{ and }
    T^* = \xi,
  \end{align*}
  where $\xi \in \real_{>0}$, $\subscr{x}{n}$, and $\lambda$ satisfy
  \begin{align}
    \label{eq:stateCostateDefender}
    \subscr{\dot x}{n} &= \xi  (A\subscr{x}{n} 
                         + \frac{ \kappa^* ~\subscr{u}{max} }{\| \subscr{p}{n} \|} B \subscr{p}{n} ) 
                         \nonumber,\\
    -  \dot \lambda &= \xi  (A^\tsp \lambda 
                      + \frac{ \kappa^* ~\subscr{u}{max} }{\| \subscr{p}{n} \|^3}  \Phi(\subscr{p}{n})
                      B^\tsp \lambda ),		 \nonumber	\\
    \dot \xi &= 0,
  \end{align}
for all $t \in [0,1]$, with boundary conditions
  \begin{align*}
    & \subscr{x}{n}(0) = x_0, \\
    & \lambda(1) = 
    2[(p(T) - \subscr{p}{n}(T))^\tsp 0_2^\tsp]^\tsp,\\
    & \lambda(0) (A x_0 +  B \kappa^*(0) \frac{p_\text{n}(0)}{\|p_\text{n} (0)\|} u_\text{max} ) = -1,
  \end{align*}
  and
  \begin{align*}
    \Phi ( \subscr{p}{n}) = 
                \left( \Vert \subscr{p}{n} \Vert^2 I_2 - 
                2\subscr{p}{n} \subscr{p}{n}^\tsp \right)
      \begin{bmatrix}
      I_2 & \zero_2\\
      \zero_2 & \zero_2
    \end{bmatrix}. 
  \end{align*}
\end{thm}
\begin{pf}
To determine the unknown final time $T$ 
we employ a technique similar to \cite{GMA-WCC:74} and 
let $t = \xi \tau$, where $\xi \in \real_{>0}$ is  a constant unknown 
parameter, and $\tau$ is the new temporal variable, with 
$0 \leq \tau \leq 1$.
We then use The Pontryagin's Maximum Principle \cite{IMG-RAS:00} 
to derive the optimality conditions for the optimization problem 
\eqref{opt:secureTrajectoryPlanning}, and 
consider the Hamiltonian
\begin{align*}
\mc H(\subscr{x}{n}, \kappa, \lambda) =  1 + 
\lambda^\tsp \xi (A\subscr{x}{n} + B \subscr{u}{n} ),
\end{align*}
where $\lambda$ is the vector function of system 
costates. 
By application of the Maximum Principle \cite{IMG-RAS:00}, the 
optimal control input and corresponding trajectory satisfy the following 
optimality conditions 
\begin{align*}
\subscr{ \dot x}{n} &= \frac{\partial \mc H}{\partial \lambda}
\Rightarrow  
\;\; \subscr{\dot x}{n} = \xi ( A\subscr{x}{n} +   B  \subscr{u}{n}  ),\\
-  \dot \lambda &= \frac{\partial \mc H}{\partial \subscr{ x}{n}}
\Rightarrow 
-  \dot \lambda = \xi ( A^\tsp \lambda + 
		 \frac{\partial  \subscr{u}{n} }{\partial  \subscr{x}{n} } B^\tsp \lambda),
\end{align*}
with boundary conditions $\subscr{x}{n}(0) \;= x_0$ and
$\lambda(1) = \frac{\partial V}{\partial x}(x(1))$,
 where we used the fact $0 \leq \tau \leq 1$.
To derive an expression for  the partial derivative of $\subscr{u}{n}$ with 
respect to $\subscr{x}{n}$ we let 
$\Pm = [e_1~e_2 ~0~0] \in \real^{2 \times 4}$ and 
rewrite \eqref{eq:secureControlInput} as
\begin{align*}
u_\text{n} = \kappa \frac{\Pm x_\text{n}}{\|\Pm x_\text{n}\|} u_\text{max},
\end{align*}
which yields
\begin{align*}
\frac{\partial  \subscr{u}{n} }{\partial  \subscr{x}{n} } = 
 \frac{\kappa \subscr{u}{max} }{\Vert \subscr{p}{n} \Vert} \Pm
- 2  \frac{\kappa \subscr{u}{max}}{\Vert \subscr{p}{n} \Vert^3} 
\Pm \subscr{x}{n} \subscr{x}{n}^\tsp \Pm^\tsp \Pm.
\end{align*} 
Hence, the expression for $\Phi ( \subscr{p}{n})$ follows by substitution.

To determine the unknown final time, we consider the additional  
differential equation $\dot \xi = 0$, and let $\xi$ be an unknown parameter.
In particular, to determine the additional boundary condition we use the fact 
that the Hamiltonian is independent of time and the final time is free.
Thus, the Hamiltonian is a first integral along optimal trajectories 
\cite{IMG-RAS:00},
i.e., $\mc H(\subscr{x}{n}, \kappa, \lambda) = \text{const.}$, with 
$\mc H(\subscr{x}{n}, \kappa, \lambda) \vert_{t=0} = 0$, which yield the 
claimed boundary conditions and the statement follows.
\QED
\end{pf}

Theorem~\ref{thm:secureTrajectoryPlanning} allows us to compute secure
control inputs by solving a two-point boundary value
problem \cite{HBK:18}. 
We propose in Fig.~\ref{fig:secureTrajectoryVsShortestPath} (see 
also corresponding data in 
Table~\ref{tab:secureTrajectoryVsShortestPath}) a comparison between a 
secure trajectory and a minimum time trajectory (non secure).
In particular, a minimum-time trajectory is obtained by numerically solving 
the following optimization problem
\begin{align*}
  & \underset{\subscr{u}{n},T}{\text{min}} 
  & & T ~+~ 
      \Vert \subscr{p}{n}(T)-\subscr{p}{F}\Vert, \nonumber \\ 
  & \text{subject to} & & \subscr{\dot x}{n} =A \subscr{x}{n} + 
  		B \subscr{u}{n}, \\
  		  &&& \subscr{p}{n}(0) = \subscr{p}{I}, \nonumber 
  \subscr{v}{n}(0) = \subscr{v}{I},\nonumber \\
\end{align*}
We note that secure trajectory require longer control horizons but, 
differently from minimum-time trajectories, prevent the existence of 
undetectable attacks.

\begin{figure}[t]
  \centering 
\subfigure[]{\includegraphics[width=\columnwidth]{%
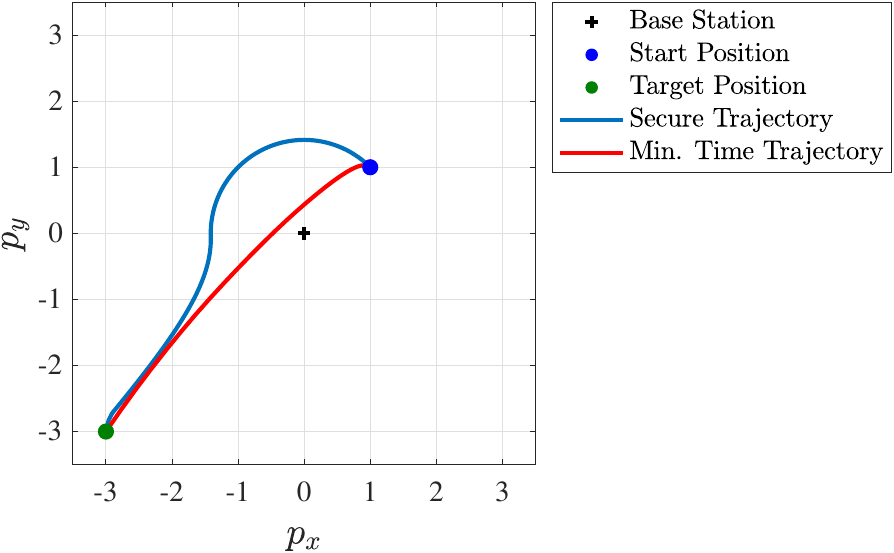}}
\subfigure[]{\includegraphics[width=.48\columnwidth]{%
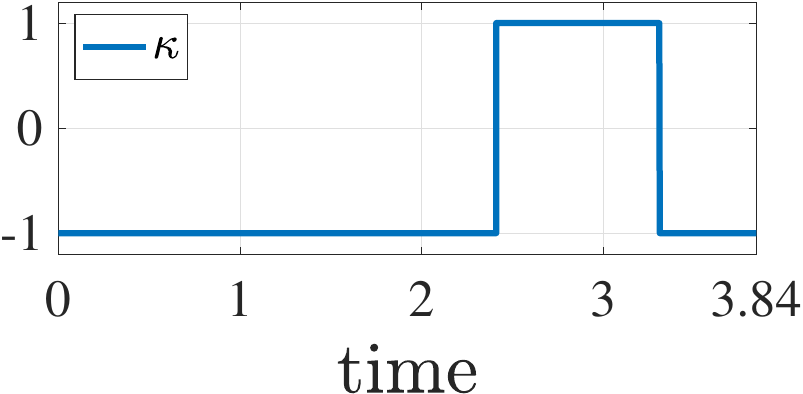}}
\subfigure[]{\includegraphics[width=.48\columnwidth]{%
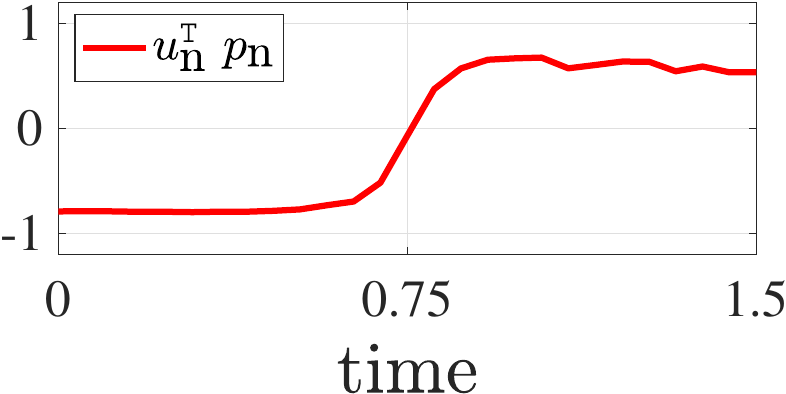}}
\label{fig: 4 int netw a} 
\caption[]{(a) Secure (blue) and minimum time (red) trajectories 
	between
      $\subscr{p}{I}=(1,1)$ and $\subscr{p}{F}=(-3,-3)$, with 
      $v(0) = [-1, 1]^\tsp$ and $\subscr{u}{max}= 1$. (b) and (c) 
      Corresponding control inputs. 
      Note that for minimum time trajectories, vectors $\subscr{u}{n}$ and 
      $\subscr{p}{n}$ are not aligned at all times.}
  \label{fig:secureTrajectoryVsShortestPath}
\end{figure}

\begin{table}[h]
\centering
\begin{tabular}{l l l}
\hline
  & \textbf{(Time)} & \textbf{(Attack Deviation)}\\
   & $T$ & $\underset{w}{\text{max}}  \Vert p(T)-p_\text{n}(T) \Vert$\\
\hline
Min. time trajectory 							& $1.5$ 		& $8.49$\\
Secure trajectory 					& $3.84$ 		& $0$\\
\hline
\end{tabular}
\caption{Details for minimum-time and secure trajectories in
Fig.~\ref{fig:secureTrajectoryVsShortestPath}.}
\label{tab:secureTrajectoryVsShortestPath}
\end{table}

\section{Undetectable attacks and secure trajectories for robots with
  unicycle dynamics}\label{sec: unicycle}
The goal of this section is to characterize undetectable attacks and
secure trajectories for robots with unicycle dynamics, so as to
illustrate that the methods developed in this paper are applicable to
a broad and general class of robot models. A robot with unicycle
dynamics has one steerable drive wheel \cite{IF-RL-MWS:00}, and
is modeled through the dynamical equations
\begin{align*}
  \dot p_\text{n}^x &= \subscr{\nu}{n} \cos (\subscr{\theta}{n}), \;
  \dot p_\text{n}^y = \subscr{\nu}{n} \sin (\subscr{\theta}{n}),\;
  \subscr{\dot \theta}{n} 	= \subscr{\omega}{n},
\end{align*}
where
$\map{\subscr{p}{n} = [p_\text{n}^x~
  p_\text{n}^y]}{\realpos}{\real^2}$ denotes the robot position,
$\map{\theta_\text{n}}{\realpos}{[0, 2 \pi)}$ denotes the steering
angle, $\map{\subscr{\nu}{n}}{\realpos}{\realpos}$ and
$\map{\subscr{\omega}{n}}{\realpos}{\real}$ denote the wheel velocity
and steering control, respectively.  
We assume that $\subscr{\nu}{n}$ is differentiable,
$\subscr{\omega}{n}$ is piecewise continuous, and
that $\subscr{\nu}{n} \leq \subscr{\nu}{max}$ and
$\vert \subscr{\omega}{n} \vert \leq \subscr{\omega}{max}$, where
$\subscr{\nu}{max} \in \real_{>0}$ and
$\subscr{\omega}{max} \in \real_{>0}$. Similarly to \eqref{eq:
  dynamics with attack}, we model the unicycle dynamics in the
presence of attacks as
\begin{align*}
\dot p^x 				&= \nu \cos (\theta), \;
\dot p^y 				= \nu \sin (\theta), \;
\dot \theta 			= \omega,
\end{align*}
where $\map{p = [p^x~ p^y]}{\realpos}{\real^2}$ and
$\map{\theta}{\realpos}{[0, 2 \pi)}$ denote, respectively, the
position and steering angle of the robot under attack, while
$\map{\nu}{\realpos}{\realpos}$ and $\map{\omega}{\realpos}{\real}$
denote the attacked wheel velocity and steering control. As described
in Section \ref{sec:problemSetup}, we assume that the robot is
equipped with a GNSS and an RSSI sensor, whose measurements are as in
\eqref{eq:nominal_readings} in the absence of attacks, and as in
\eqref{eq:attack_readings} in the presence of attacks.

Using the notions in Definition~\ref{def: undetectable through y1 and
  y2}, we next characterize undetectable attacks against robots with
unicycle dynamics. 
In the remainder of this section, we let $\operatorname{angle}(v,w)$ 
denote the angle between the vectors $v$ and~$w$, that is 
\begin{align*}
\operatorname{angle}(v,w) = 
\arccos \left( \frac{v^\tsp w}{\Vert v \Vert \Vert w \Vert} \right) ,
\end{align*}
with $\operatorname{angle}(0,w)=\operatorname{angle}(v,0)=0$.

\begin{thm}{\bf \textit{(Undetectable attacks for unicycle dynamics)}}
\label{thm:undetectableAttacksUnicycle}
  Let
  $\subscr{\phi}{n} = \operatorname{angle} ( \subscr{p}{n} , 
  \subscr{\dot p}{n})$ and $\phi = \operatorname{angle} (p , \dot p)$. 
  The attack $(\nu, \omega, \supscr{u}{GNSS})$ is undetectable if and only
  if
  \begin{align}
    \label{eq:undetectabilityUnicycle}
    \nu \cos (\phi) = \subscr{\nu}{n} \cos ( \subscr{\phi}{n}),
    \text{ and }   \supscr{u}{GNSS} = \subscr{p}{n}-p .
  \end{align}
\end{thm}
\begin{pf}
  The proof follows by extending the proof of Theorem~\ref{thm: main
    zero dynamics} to the considered unycicle dynamics.  \QED
\end{pf}

\begin{figure}[tb]
  \centering
  \includegraphics[width=.6\columnwidth]{%
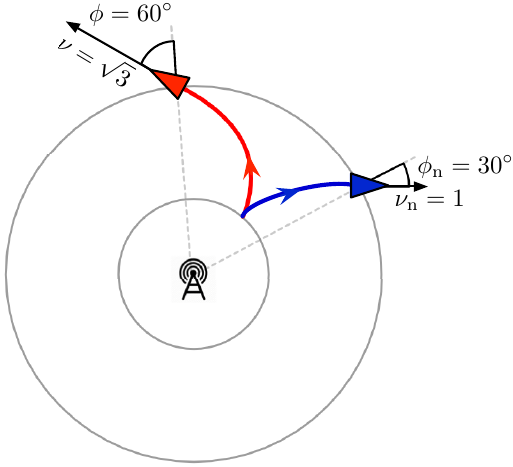}
  \caption{Nominal (blue) and undetectable attack (red) trajectories. 
  As discussed in Theorem
    \ref{thm:undetectableAttacksUnicycle}, the illustrated vectors
    satisfy
    $\nu \cos (\phi) = \subscr{\nu}{n} \cos ( \subscr{\phi}{n})$ to
    guarantee undetectability.}
    \label{fig:undetectableTrajectoriesUnicycle}
\end{figure}

An example where the condition in 
Theorem~\ref{thm:undetectableAttacksUnicycle} holds is
illustrated in Fig.~\ref{fig:undetectableTrajectoriesUnicycle}. 
Next, we provide a characterization of secure control inputs for unicycle 
dynamics.

\begin{thm}{\bf \textit{(Secure control inputs for 
      unicycle dynamics)}}\label{thm:secureControlInputUnicycle}
Let $\sbs{\phi}{n} = \operatorname{angle} ( \sbs{p}{n} ,  \sbs{\dot p}{n})$.
The control input $(\subscr{\nu}{n}, \subscr{\omega}{n})$ is secure if and 
only if the following conditions  hold simultaneously:
  \begin{enumerate}
  \item[(1)] \blue{$\subscr{\phi}{n}(0) \in \{0, \pi \}$}, and $p_\text{n} \neq 0$ at all 
  					times, 
  \item[(2)] $\nu_\text{n} = \nu_\text{max}$ and
    $\omega_\text{n} = 0$ at all times.
  \end{enumerate}
\end{thm}

\begin{pf}
(\textit{If}) 
We assume (1)-(2) and show that either (C1) or (C2) is satisfied. 
We distinguish among two cases. \\
\textit{(Case 1)} The attack $(\nu, \omega, \supscr{u}{GNSS})$ 
does not satisfy \eqref{eq:undetectabilityUnicycle}. 
Then, (C2)  immediately follows. \\
\textit{(Case 2)} The attack $(\nu, \omega, \supscr{u}{GNSS})$ satisfies
\eqref{eq:undetectabilityUnicycle}. 
We first focus on the time instant $\tau = 0$, and take the absolute value on 
both sides of the undetectability condition \eqref{eq:undetectabilityUnicycle}  
to obtain the identity
\begin{align}
\label{eq:auxThm14}
\vert \subscr{\nu}{n} (0) \cos ( \subscr{\phi}{n}(0))\vert = 
\vert \nu(0) \cos (\phi(0)) \vert.
\end{align}
By substituting conditions (1)-(2) into the left-hand side of 
\eqref{eq:auxThm14} we obtain
\begin{align*}
\vert \sbs{\nu}{n}(0) \cos (\sbs{\phi}{n}(0)) \vert = \sbs{\nu}{max}.
\end{align*}
On the other hand, by applying the Cauchy-Schwartz inequality to the 
right-hand side of \eqref{eq:auxThm14} we have 
\begin{align*}
\vert \nu(0) \cos (\phi(0)) \vert \leq 
\vert \nu(0) \vert \vert \cos (\phi (0)) \vert \leq 
\sbs{\nu}{max}.
\end{align*}
By application of \eqref{eq:auxThm14}, exact equality must hold in the above 
bound, and thus we necessarily have 
$\nu(0) = \sbs{\nu}{max}$ and  $\vert \cos (\phi (0)) \vert =1$.
Finally, we use the fact that $\phi$ is a differentiable function of time 
to obtain $\phi(0) = \sbs{\phi}{n}(0)$.
To conclude the proof, we use the fact that $\sbs{\omega}{n}=0$ at all 
times, which implies $\sbs{\phi}{n}(\tau) = \sbs{\phi}{n}(0)$ for all $\tau$
(a formal proof of this fact can be done by leveraging the change of 
variables \eqref{eq:changeOfVariables} below), and iterate the above 
reasoning for all $\tau \in [0,T)$ to obtain 
$\nu = \subscr{\nu}{n}$ and $\phi= \subscr{\phi}{n}$ at all times, from 
which condition (C1) follows.

(\textit{Only if})
We now show that (C1)-(C2) imply (1)-(2)  or, 
equivalently, if (1)-(2) do not hold, then (C1)-(C2) are not verified.
We distinguish among four cases.\\
\textit{(Case 1)}
The robot initial conditions are such that 
$\subscr{\theta}{n}(0) \neq \{ 0, \pi \}$, and the nominal control inputs satisfy 
$\subscr{\nu}{n} = \subscr{\nu}{max}$ and $ \subscr{\omega}{n}(t) = 0$ 
at all times.
We consider the attack $(\nu, \omega)$, with 
\begin{align*}
\dot \nu = \frac{1}{\cos \phi} 
\left( 
\subscr{\nu}{n}  \subscr{\dot \phi}{n} \sin \subscr{\phi}{n}   - 
\nu \dot \phi \sin \dot \phi 
\right),
 \text{ and } 
\omega \neq 0,
\end{align*}
and show that such attack is undetectable and violates (C1)-(C2).
To prove undetectability, we use the fact that 
$\subscr{x}{n}(0) = x(0)$,
and equivalently show the identity between the time derivatives of 
\eqref{eq:undetectabilityUnicycle}, which yields
\begin{align*}
- \subscr{\nu}{n}  \subscr{\dot \phi}{n} \sin \subscr{\phi}{n} = 
\dot \nu \cos \phi - \nu \dot \phi \sin \phi,
\end{align*}
where we used the relationship $\subscr{\dot \nu}{n} = 0$ at all times.
As a result, undetectability of the given attack $(\nu, \omega)$
follows by substitution. 
To conclude, we observe that $(\nu, \omega)$ is undetectable and 
violates (C1) and (C2).\\
\textit{(Case 2)}
There exists $\bar t \in [0,T]$ such that $\| \subscr{p}{n}( \bar t) \| = 0$.
It follows from \eqref{eq:undetectabilityUnicycle} that, when 
$\subscr{p}{n}(\bar t) =0$,
$\operatorname{angle}(0,  \subscr{\dot p}{n}(\bar t) ) = 0$.
As a result, attack inputs are unconstrained at time $\bar t$, and any 
$\omega$ such that $\omega(t) = \subscr{\omega}{n}(t)$ for all 
$t \in [0,\bar t)$, and $\omega (\bar t) \neq \subscr{\omega}{n}(\bar t)$, is 
undetectable and violates (C1)-(C2).\\
\textit{(Case 3)}
Nominal control inputs satisfy $\subscr{\nu}{n} = \subscr{\nu}{max}$ at all 
times, $ \subscr{\omega}{n}(t) = 0$ for all $t \in [0, \bar t)$, and
$ \subscr{\omega}{n}( \bar t) \neq 0$, $\bar t \in [0,T]$.
We perform a change of variables  \cite{BS-LS-LV-GO:10}
\begin{align}
\label{eq:changeOfVariables}
\rho 				&= \sqrt{p_x^2 + p_y^2}, \nonumber\\
\phi  			&= \operatorname{Atan2}(p_y, p_x) - \theta + \pi, \nonumber\\
\delta 			&= \phi + \theta,
\end{align}
which leads to the following dynamical equation 
$\dot \phi = \frac{\sin \phi}{\Vert p \Vert} - \omega$ 
that relates $\phi$ to the control inputs $\nu$ and $\omega$
(see \cite{BS-LS-LV-GO:10}).
We then let $\nu = \subscr{\nu}{n}$ at all times, and
\begin{align*}
\omega = 
\begin{cases}
\subscr{\omega}{n}, 			& \text{if } t \in [0,\bar t),\\
- \subscr{\omega}{n}, 			& \text{if } t \in [\bar t, T],\\
\end{cases}
\end{align*}
from which we obtain
\begin{align*}
\dot \phi = 
\frac{\sin \phi}{\Vert p \Vert} - \omega =
- \frac{\sin\subscr{\phi}{n} }{\Vert \subscr{p}{n}\Vert} + \subscr{\omega}{n}  =
- \subscr{\dot \phi}{n},
\end{align*}
for all $t \in [\bar t, T]$.
By combining the above relationship with 
$\subscr{\phi}{n}(\bar t) = \phi (\bar t)$ we obtain  
$\phi ( t ) = - \subscr{\phi}{n}( t)$ for all $t \in [\bar t, T]$.
To conclude, we note that the given choice of $(\nu, \omega)$ leads to
an undetectable attack that satisfies \eqref{eq:undetectabilityUnicycle} and 
that violates (C1)-(C2).\\
\textit{(Case 4)}
Nominal control inputs satisfy $\subscr{\omega}{n} = 0$ at all times, 
$ \subscr{\nu}{n}(t) = \subscr{\nu}{max}$ for all $t \in [0, \bar t)$, and
$ \subscr{\nu}{n}( \bar t) < \subscr{\nu}{max}$, $\bar t \in [0,T]$.
Under these assumptions, we consider the attack $(\nu, \omega)$ with
$\nu(t) = \subscr{\nu}{n}(t)$ and $\omega(t) = \subscr{\omega}{n}(t)$ 
for all $t \in [0, \bar t)$, and
\begin{align*}
\dot \nu = \frac{1}{\cos \phi} \left( \subscr{\dot \nu}{n} +
\nu \dot \phi \sin \phi \right),
 \text{ and } 
\omega \neq 0,
\end{align*}
for all $t \in [\bar t, T]$.
To prove undetectability of the considered attack, we observe that 
$\subscr{x}{n}(\bar t) = x(\bar t)$,
and equivalently show the identity between the time derivatives 
of \eqref{eq:undetectabilityUnicycle} for all $t \in (\bar t, T]$, which reads
\begin{align*}
\subscr{\dot \nu}{n}  = \dot \nu \cos \phi - \nu \dot \phi \sin \phi,
\end{align*}
where we used the relationships $\cos \subscr{\phi}{n} = 1$ and 
$\sin \subscr{\phi}{n}=0$ at all times.
As a result, undetectability of the given attack follows by substitution. 
To conclude, we note that $(\nu, \omega)$ is undetectable and 
violates (C1) and (C2). 
\QED
\end{pf}
From Theorem~\ref{thm:secureControlInputUnicycle}, a secure trajectory 
exists only if initial position, final position, and the origin are collinear.
We conclude by observing that this aspect is consistent with the similar 
conclusions previously drawn in Theorem~\ref{thm:secureControlInput} and 
Lemma~\ref{lem:accessibleConfigurations}.

\begin{rem}{\bf \textit{(Generality of our methods)}}
  The approach presented in this work for the characterization of
  undetectable attacks and the resulting effects on the robot trajectories 
  can be generalized to a
  wider class of nonlinear systems and attacks. The approach consists
  of three main steps, namely, the characterization of undetectable
  trajectories by studying the Lie derivatives of the measurement
  equations, the characterization of undetectable inputs and secure
  trajectories by solving a set of nonlinear algebraic equations akin
  to \eqref{eq: implicit undetectability condition} and
  \eqref{eq:undetectabilityUnicycle}, and the study of the submanifold
  of the state space that can be reached by undetectable
  attacks. While systematic tools may exist to solve the first two
  steps for a broad class of dynamics, the problem of nonlinear
  constrained controllability in the third step, which is solved in
  this paper via numerical optimization, requires the development of
  novel control theories and tools.
\QEDB
\end{rem}

\section{Conclusions}\label{sec:conclusions}
In this paper we introduce and study the problem of secure trajectory
planning, that is, the design of trajectories to guarantee the navigation
between two desired configurations despite the presence of attackers.
We focus on the case where
the robot has a GNSS sensor and a RSSI sensor, and provide an explicit
characterization of secure trajectories, undetectable attacks, and
their effects on the nominal trajectory. Further, we provide
numerical algorithms to determine secure trajectories and optimal
attacks, and we illustrate our findings through a set of examples. 
To the best of our knowledge, this work constitutes a first step 
towards understanding the fundamental limitations of attack detection 
algorithms for systems with nonlinear dynamics.
Several aspects are left as the subject of future investigation, including a
characterization of the set of configurations reachable by
undetectable attacks and secure trajectories, the extension of the
methods to different classes of sensors and attacks, the development 
of robust control mechanisms to operate the system despite the presence 
of attacks, and the study of
secure trajectories and undetectable attacks in the presence of
sensing and actuation noise.



\bibliographystyle{elsarticle-num}
\bibliography{alias,FP,Main,New}

\begin{thebibliography}{10}
\expandafter\ifx\csname url\endcsname\relax
  \def\url#1{\texttt{#1}}\fi
\expandafter\ifx\csname urlprefix\endcsname\relax\def\urlprefix{URL }\fi
\expandafter\ifx\csname href\endcsname\relax
  \def\href#1#2{#2} \def\path#1{#1}\fi

\bibitem{YZL-AD-FS-IM-MDDB:19}
Y.~Z. Lun, A.~D'Innocenzo, F.~Smarra, I.~Malavolta, M.~D.~D. Benedetto, State
  of the art of cyber-physical systems security: An automatic control
  perspective, The Journal of Systems and Software 149~(2019) (2019) 174--216.

\bibitem{FP-FD-FB:10y}
F.~Pasqualetti, F.~D{\"o}rfler, F.~Bullo, Attack detection and identification
  in cyber-physical systems, IEEE Transactions on Automatic Control 58~(11)
  (2013) 2715--2729.

\bibitem{CZB-FP-VG:16}
C.-Z. Bai, F.~Pasqualetti, V.~Gupta, Data-injection attacks in stochastic
  control systems: Detectability and performance tradeoffs, Automatica 82
  (2017) 251--260.

\bibitem{YM-BS:10a}
Y.~Mo, B.~Sinopoli, Secure control against replay attacks, in: Allerton Conf.\
  on Communications, Control and Computing, Monticello, IL, USA, 2010, pp.
  911--918.

\bibitem{FH-PT-SD:11}
F.~Hamza, P.~Tabuada, S.~Diggavi, Secure state-estimation for dynamical systems
  under active adversaries, in: Allerton Conf.\ on Communications, Control and
  Computing, Monticello, IL, USA, 2011, pp. 337--344.

\bibitem{JPH-SDB:19}
J.~P. Hespanha, S.~D. Bopardikar, Output-feedback linear quadratic robust
  control under actuation and deception attacks, in: {A}merican {C}ontrol
  {C}onference, Philadelphia, PA, USA, 2019, pp. 489--496.

\bibitem{YS-PN-NB-AS-SA-PT-15}
Y.~Shoukry, P.~Nuzzo, N.~Bezzo, A.~L. Sangiovanni-Vincentelli, S.~A. Seshia,
  P.~Tabuada, Secure state reconstruction in differentially flat systems under
  sensor attacks using satisfiability modulo theory solving, in: {IEEE} Conf.\
  on Decision and Control, 2015, pp. 3804--3809.

\bibitem{QU-DF-YHC-CJT:18}
Q.~Hu, D.~Fooladivanda, Y.~H. Chang, C.~J. Tomlin, Secure state estimation and
  control for cyber security of the nonlinear power systems, IEEE Transactions
  on Automatic Control 5~(3) (2018) 1310--1321.

\bibitem{JK-CL-HS-YE-JHS:18}
J.~Kim, C.~Lee, H.~Shim, Y.~Eun, J.~H. Seo, Detection of sensor attack and
  resilient state estimation for uniformly observable nonlinear systems having
  redundant sensors, IEEE Transactions on Automatic Control 64~(3) (2019)
  1162--1169.

\bibitem{AB-AJ-VD-JN-GL:12}
A.~Broumandan, A.~Jafarnia-Jahromi, V.~Dehghanian, J.~Nielsen, G.~Lachapelle,
  {GNSS} spoofing detection in handheld receivers based on signal spatial
  correlation, in: ION Position, Location and Navigation Symposium, 2012, pp.
  479--487.

\bibitem{XJ-JZ-BJH-JJM-ADD:13}
X.~Jiang, J.~Zhang, B.~J. Harding, J.~J. Makela, A.~D.
  Dom{\'\i}nguez-Garc{\'\i}a, Spoofing {GPS} receiver clock offset of phasor
  measurement units, IEEE Transactions on Power Systems 28~(3) (2013)
  3253--3262.

\bibitem{PFS-BNA-KCS-RJH:14}
P.~F. Swaszek, S.~A. Pratz, B.~N. Arocho, K.~C. Seals, R.~J. Hartnett, {GNSS}
  spoof detection using shipboard {IMU} measurements, in: International
  Technical Meeting of The Satellite Division of the Institute of Navigation,
  Tampa, FL, 2014, pp. 745--758.

\bibitem{QZ-SH-FL-MC:16}
Q.~Zou, S.~Huang, F.~Lin, M.~Cong, Detection of {GPS} spoofing based on {UAV}
  model estimation, in: Annual Conference of the IEEE Industrial Electronics
  Society, 2016, pp. 6097--6102.

\bibitem{DSR-PFS-KCS:15}
D.~S. Radin, P.~F. Swaszek, K.~C. Seals, {GNSS} spoof detection based upon
  pseudoranges from multiple receivers, in: International Technical Meeting of
  The Institute of Navigation, Dana Point, CA, 2015, pp. 657--671.

\bibitem{LMP-BWO-SPP-JAB-KDW-TEH:14}
M.~L. Psiaki, B.~W. O'Hanlon, S.~P. Powell, J.~A. Bhatti, K.~D. Wesson, T.~E.
  Humphreys, A.~Schofield, {GNSS} spoofing detection using two-antenna
  differential carrier phase, in: International Technical Meeting of the
  Satellite Division of The Institute of Navigation, Tampa, FL, 2014, pp.
  2776--2800.

\bibitem{MLP-BWO-JAB-DPS-TEH:13}
M.~L. Psiaki, B.~W. O'Hanlon, J.~A. Bhatti, D.~P. Shepard, T.~E. Humphreys,
  {GPS} spoofing detection via dual-receiver correlation of military signals,
  IEEE Transactions on Aerospace and Electronic System 49~(4) (2013)
  2250--2267.

\bibitem{PYM-TEH-NML:09}
P.~Y. Montgomery, T.~E. Humphreys, B.~M. Ledvina, Receiver-autonomous spoofing
  detection: Experimental results of a multi-antenna receiver defense against a
  portable civil {GPS} spoofer, in: ION International Technical Meeting, 2009,
  pp. 124--130.

\bibitem{MJ-RD:03}
M.~Jun, R.~{D'Andrea}, Path planning for unmanned aerial vehicles in uncertain
  and adversarial environments, in: Cooperative control: models, applications
  and algorithms, Springer, USA, 2003, pp. 95--110.

\bibitem{DPS-TEH-AAF:12}
D.~P. Shepard, T.~E. Humphreys, A.~A. Fansler, Evaluation of the vulnerability
  of phasor measurement units to {GPS} spoofing attacks, International Journal
  of Critical Infrastructure Protection 5~(3-4) (2012) 146--153.

\bibitem{AJK-DPS-JAB-TEH:14}
A.~J. Kerns, D.~P. Shepard, J.~A. Bhatti, T.~E. Humphreys, Unmanned aircraft
  capture and control via {GPS} spoofing, Journal of Field Robotics 31~(4)
  (2014) 617--636.

\bibitem{KH-CS-13}
K.~Hartmann, C.~Steup, The vulnerability of {UAVs} to cyber attacks - an
  approach to the risk assessment, in: International Conference on Cyber
  Conflict, 2013, pp. 1--23.

\bibitem{SML:06}
S.~M. LaValle, Planning Algorithms, Cambridge University Press, 2006, available
  at \texttt{http://planning.cs.uiuc.edu}.

\bibitem{TB-GJO:99}
T.~Ba{\c s}ar, G.~J. Olsder, Dynamic Noncooperative Game Theory, 2nd Edition,
  SIAM, 1999.

\bibitem{GB-FP:19}
G.~Bianchin, Y.-C. Liu, F.~Pasqualetti, Secure navigation of robots in
  adversarial environments, IEEE Control Systems Letters 4~(1) (2019) 1--6.

\bibitem{IMG-RAS:00}
I.~M. Gelfand, R.~A. Silverman, et~al., Calculus of variations, Courier
  Corporation, 2000.

\bibitem{RFH-SPS-RGV:95}
R.~F. Hartl, S.~P. Sethi, R.~G. Vickson, A survey of the maximum principles for
  optimal control problems with state constraints, SIAM Review 37~(2) (1995)
  181--218.

\bibitem{HBK:18}
H.~B. Keller, Numerical methods for two-point boundary-value problems, Courier
  Dover Publications, 2018.

\bibitem{GMA-WCC:74}
G.~M. Aly, W.~C. Chan, Numerical computation of optimal control problems with
  unknown final time, Journal of Mathematical Analysis and Applications 45~(2)
  (1974) 274--284.

\bibitem{IF-RL-MWS:00}
I.~Fantoni, R.~Lozano, M.~W. Spong, Energy based control of the pendubot, IEEE
  Transactions on Automatic Control 45~(4) (2000) 725--729.

\bibitem{BS-LS-LV-GO:10}
B.~Siciliano, L.~Sciavicco, L.~Villani, G.~Oriolo, Robotics: modelling,
  planning and control, Springer Science \& Business Media, 2010.

\end{thebibliography}

\end{document}